\documentclass[11pt,reqno]{amsart}
\usepackage{graphicx}

\usepackage{amssymb}
\usepackage{cite}
\usepackage{amsmath}
\usepackage{latexsym}
\usepackage{amscd}
\usepackage{amsthm}
\usepackage{mathrsfs}
\usepackage{url}

\makeatletter
\@namedef{subjclassname@2020}{%
  \textup{2020} Mathematics Subject Classification}
\makeatother

\usepackage[backref=page]{hyperref}
\numberwithin{equation}{section}

\newtheorem{thm}{Theorem}

\newtheorem{lem}[thm]{Lemma}
\newtheorem{prop}[thm]{Proposition}

\newtheorem{exam}{Example}[section]

\theoremstyle{definition}

\theoremstyle{remark}
\newtheorem{rem}[thm]{Remark}

\def\C{\mathbb C}
\def\R{\mathbb R}

\def\N{\mathbb N}

\def\SS{\mathbb S}

\def\pt{\partial}

\allowdisplaybreaks

\begin{document}
\title[Sharp estimates of Steklov eigenvalues]{Sharp estimates of Steklov eigenvalues for differential forms on warped product manifolds}
\author[Q.~Liu]{Qinyong~Liu}
\address{School of Mathematics, Sichuan University, Chengdu 610065, Sichuan,  P.~R.~China}
\email{\href{mailto:1197078867@qq.com}{1197078867@qq.com}}
\author[C.~Xiong]{Changwei~Xiong}
\address{School of Mathematics, Sichuan University, Chengdu 610065, Sichuan,  P.~R.~China}
\email{\href{mailto:changwei.xiong@scu.edu.cn}{changwei.xiong@scu.edu.cn}}
\date{\today}

\maketitle

\begin{abstract}
We consider the Steklov eigenvalue problem on differential forms of degree $p$ introduced by Karpukhin (2019) (and by Belishev and Sharafutdinov (2008)) on $n$-dimensional warped product manifolds with nonnegative Ricci curvature and strictly convex boundary. We extend Chakradhar's lower bound for the spectrum of Escobar type (2025) to the full range $p\in [1,n-2]$ and improve on his result for $p\in [1,(n-2)/2]$. We also obtain sharp estimates for the gaps and ratios of consecutive eigenvalues for $p\in [1,(n-2)/2]$.
\end{abstract}

\section{Introduction}\label{sec1}

In this paper we study sharp estimates of Steklov eigenvalues for differential forms on warped product manifolds. The classical Steklov eigenvalue problem for functions (i.e., $0$-forms) was introduced by Steklov \cite{Ste02} around 1900 (cf. \cite{KKK14}) and has received lots of attention since then. See \cite{GP17, CGGS24} for surveys on Steklov eigenvalue problems. In contrast, the studies on Steklov problems for differential forms are relatively few and in some sense harder than those for functions. Here we add some new results for differential forms which we hope would motivate further investigation on the topic.

In the literature, roughly speaking, there are three types of Steklov eigenvalue problems (or the Dirichlet-to-Neumann maps) for differential forms. One is introduced and studied by Belishev and Sharafutdinov \cite{BS08} and Karpukhin \cite{Kar19}. A second one, the complete Dirichlet-to-Neumann map, is by Joshi and Lionheart \cite{JL05} and Sharafutdinov and Shonkwiler \cite{SS13}. (In fact, this one determines the previous one.) Another one is by Raulot and Savo \cite{RS12} (cf. Carron's relative to absolute operator \cite{Car02}). In this paper we are mainly concerned with the one in \cite{BS08,Kar19}. Precisely, the problem may be described as follows. Let $(M^{n},g)$ be an $n$-dimensional ($n\geq 3$) smooth compact connected oriented Riemannian manifold with boundary $\Sigma=\pt M$. Given an integer $0\leq p\leq n-1$ and a differential form $\phi \in \Omega^p(\Sigma)$, there exists a solution $\varphi\in \Omega^p(M)$ satisfying
\begin{align}
\begin{cases}
\Delta \varphi=0,\ \delta \varphi  =0\text{ on }M,\\
J^*\varphi =\phi\text{ on }\Sigma,
\end{cases}
\end{align}
where $J:\Sigma\to M$ is the inclusion mapping. Then the Dirichlet-to-Neumann operator $\Lambda^{[p]}$ is defined by
\begin{equation*}
\Lambda^{[p]}(\phi):=i_\nu d\varphi\in \Omega^p(\Sigma),
\end{equation*}
where $\nu$ denotes the unit outward normal vector field along the boundary $\Sigma$. In order to study the spectral properties of $\Lambda^{[p]}$, it is natural to restrict $\Lambda^{[p]}$ to the subspace of co-closed $p$-forms in $\Omega^p(\Sigma)$ (see \cite{Kar19}). In that case the restricted operator admits a sequence of eigenvalues
\begin{align*}
0\leq \sigma_1^{[p]}(\Sigma)\leq \sigma_2^{[p]}(\Sigma)\leq \dots \nearrow \infty,
\end{align*}
with the account of multiplicities. Throughout the paper we also use $\sigma^{[p]}_{(m)}$ ($m\in \N$) to denote the eigenvalue without counting multiplicities.

Obtaining estimates as sharp as possible on the spectrum $\sigma_m^{[p]}(\Sigma)$ under various geometric conditions would help us understand the relationship of the geometry and topology of the underlying manifolds, as for classical eigenvalue problems. In this aspect, Karpukhin \cite{Kar19} proved a Hersch--Payne--Schiffer type inequality for $\Lambda^{[p]}$ on differential forms. Girouard, Karpukhin, Levitin, and Polterovich \cite{GKLP22} obtained an estimate comparing the eigenvalues of $\Lambda^{[p]}$ on co-closed forms with the corresponding eigenvalues of the Hodge Laplacian on the boundary. For Riemannian manifolds with certain curvature conditions, the second-named author \cite{Xio26} derived a sharp lower bound for the first positive eigenvalue of $\Lambda^{[p]}$ on differential forms, and gave another comparison result between the spectrum of $\Lambda^{[p]}$ and the spectrum of the Hodge Laplacian on the boundary. While for warped product manifolds with nonnegative Ricci curvature and strictly convex boundary, Chakradhar \cite{Cha25} deduced a sharp lower bound for the spectrum of $\Lambda^{[p]}$ when $1\leq p\leq (n-1)/2$.

In this paper we continue the study of Chakradhar's on warped product manifolds. Although the class of warped product manifolds is of high symmetry, it includes important examples of space forms, Schwarzschild manifolds, anti-de Sitter-Schwarzschild manifolds, etc. of physical significance. Moreover, the results regarding eigenvalue problems on warped product manifolds sometimes may provide useful hints for studies on general Riemannian manifolds. So in this work we restrict ourselves to warped product manifolds; see Section~\ref{sec2.3} for the precise setting. First we extend Chakradhar's result \cite{Cha25} to the full range of the degree $p\in [1,n-2]$ of the differential forms.
\begin{thm}\label{thm1}
Let $M=[0,R]\times\mathbb{S}^{n-1}$ be an $n$-dimensional ($n\geq 3$) warped product manifold with the metric
\begin{equation*}
    g=dr^{2}+h^{2}(r)g_{\mathbb{S}^{n-1}},
\end{equation*}
where $h(r)$ satisfies the smoothness assumption~\eqref{eq-A}. Suppose that $M$ has $\mathrm{Ric}_g\geq 0$ and the boundary $\pt M$ is strictly convex. Then for any $1\leq p\leq n-2$, we have
    \begin{equation}
        \sigma^{[p]}_{(m)}\geq (m+p)\frac{h'(R)}{h(R)},\ \forall m\geq 1.
    \end{equation}
The equality holds for any given $m\in \N$ if and only if $h(r)=r$, i.e., $M$ is isometric to a Euclidean ball of radius $R$.
\end{thm}
\begin{rem}\label{rem2}
For $p\in [1,(n-1)/2]$, the result in Theorem~\ref{thm1} is due to Chakradhar \cite{Cha25} as mentioned above. Besides, for $p=n-1$, there is a unique eigenvalue $\sigma^{[n-1]}_{(1)}$ which satisfies $nh'(R)/h(R)\leq \sigma^{[n-1]}_{(1)}\leq n/h(R)$; see Remarks~\ref{rem9} and \ref{rem13} for details. We also note that $h'(R)/h(R)$ is the principal curvature of the boundary, and so the result in Theorem~\ref{thm1} is of Escobar type (cf. \cite{Esc97,Esc99,XX24,Xio26}).
\end{rem}

Second, for small degrees $p$ of the differential forms, we are able to get a sharper estimate than that in Theorem~\ref{thm1} (note $h'(R)\leq 1$ by Lemma~\ref{lem-h} below).
\begin{thm}\label{thm2}
Let $M=[0,R]\times\mathbb{S}^{n-1}$ be an $n$-dimensional ($n\geq 4$) warped product manifold with the metric
\begin{equation*}
    g=dr^{2}+h^{2}(r)g_{\mathbb{S}^{n-1}},
\end{equation*}
where $h(r)$ satisfies the smoothness assumption~\eqref{eq-A}. Suppose that $M$ has $\mathrm{Ric}_g\geq 0$ and the boundary $\pt M$ is strictly convex. Then for any $1\leq p\leq (n-2)/2$, we have
    \begin{equation}\label{eq-1.3}
        \sigma^{[p]}_{(m)}\geq \frac{m+p}{h(R)},\ \forall m\geq 1.
    \end{equation}
If $p< (n-2)/2$, then the equality for any given $m\in \N$ holds if and only if $h(r)=r$, i.e., $M$ is isometric to a Euclidean ball of radius $R$. If $p= (n-2)/2$ is an integer, then the equality always holds, i.e.,
    \begin{equation}\label{eq-1.4}
        \sigma^{[p]}_{(m)}=\frac{m+p}{h(R)}, \ \forall m\geq 1.
    \end{equation}
\end{thm}
\begin{rem}
When $p=(n-2)/2$ is an integer, from the proof we can see that \eqref{eq-1.4} holds without the assumptions on $\mathrm{Ric}_g$ and the convexity of $\pt M$. In fact, the result for $p=(n-2)/2$ is due to Chakradhar \cite[Theorem~1.4 (ii)]{Cha25} where he used a conformal argument.
\end{rem}
\begin{rem}
The lower bound in \eqref{eq-1.3} can be interpreted geometrically. For instance, we may think of $h(R)$ as $h(R)=(|\pt M|/|\SS^{n-1}|)^{1/(n-1)}$, a normalized area of the boundary.
\end{rem}
Last, the investigation on sharp estimates for eigenvalue gaps and ratios is also an important topic in spectral geometry. See, e.g., \cite{AC11,AB89,AB92,AB01,BL07,Kro92,Lav94,PPW56,PW60,ZY84} for some results on the classical Dirichlet and Neumann eigenvalue problems. Here for the Steklov eigenvalue problem on differential forms, we can derive sharp estimates for the gaps and ratios of consecutive eigenvalues as follows.
\begin{thm}\label{thm3}
Under the conditions of Theorem~\ref{thm2}, for any $1\leq p\leq (n-2)/2$ we have
    \begin{equation}
        \sigma^{[p]}_{(m+1)}-\sigma^{[p]}_{(m)}\geq \frac{1}{h(R)},\ \forall m\geq 1.
    \end{equation}
If $p<(n-2)/2$, then the equality for any given $m\in \N$ holds if and only if $h(r)=r$, i.e., $M$ is isometric to a Euclidean ball of radius $R$. If $p= (n-2)/2$ is an integer, then we always have
    \begin{equation}
        \sigma^{[p]}_{(m)}=\frac{m+p}{h(R)}, \ \forall m\geq 1.
    \end{equation}
\end{thm}
\begin{thm}\label{thm4}
Under the conditions of Theorem~\ref{thm2}, for any $1\leq p\leq (n-2)/2$ we have
    \begin{equation}
        \frac{\sigma^{[p]}_{(m+1)}}{\sigma^{[p]}_{(m)}}\leq\frac{m+p+1}{m+p}, \ \forall m\geq 1.
    \end{equation}
If $p<(n-2)/2$, then the equality for any given $m\in \N$ holds if and only if $h(r)=r$, i.e., $M$ is isometric to a Euclidean ball of radius $R$. If $p= (n-2)/2$ is an integer, then we always have
    \begin{equation}
        \sigma^{[p]}_{(m)}=\frac{m+p}{h(R)}, \ \forall m\geq 1.
    \end{equation}
\end{thm}
\begin{rem}
For $0$-forms or functions ($p=0$), the results as in Theorems~\ref{thm3} and \ref{thm4} were derived by the second-named author in \cite{Xio21}.
\end{rem}
As for the proofs of Theorems~\ref{thm1}, \ref{thm2}, \ref{thm3} and \ref{thm4}, we use Chakradhar's characterization \cite{Cha25} of eigenvalues and eigenforms of the operator $\Lambda^{[p]}$ on warped product manifolds. Together with the properties of the warping function $h(r)$ resulted from the curvature assumptions in the theorems, we can draw the conclusions.

Adapting the proofs of above results, we can handle the case of $\mathrm{Ric}_{g}\leq 0$. More precisely, we have the following remark.
\begin{rem}\label{rem6}
Assume the same conditions of Theorem~\ref{thm2} except $\mathrm{Ric}_{g}\geq 0$ replaced by $\mathrm{Ric}_{g}\leq 0$. Then for any $1\leq p\leq(n-2)/2$, we have
    \begin{equation*}
        \sigma^{[p]}_{(m)}\leq \frac{m+p}{h(R)},\ \forall m\geq 1,
    \end{equation*}
    \begin{equation*}
        \sigma^{[p]}_{(m+1)}-\sigma^{[p]}_{(m)}\leq\frac{1}{h(R)},\ \forall m\geq 1,
    \end{equation*}
and
    \begin{equation*}
        \frac{\sigma^{[p]}_{(m+1)}}{\sigma^{[p]}_{(m)}}\geq\frac{m+p+1}{m+p},\ \forall  m\geq 1.
    \end{equation*}
If $p<(n-2)/2$, then the equality in any of the above inequalities holds if and only if $h(r)=r$, i.e., $M$ is isometric to a Euclidean ball of radius $R$. If $p= (n-2)/2$ is an integer, then we always have
    \begin{equation*}
        \sigma^{[p]}_{(m)}=\frac{m+p}{h(R)}, \ \forall m\geq 1.
    \end{equation*}
A similar conclusion can be drawn for Theorem~\ref{thm1}.
\end{rem}

The paper is organized as follows. In Section~\ref{sec2} we present some preliminaries on differential forms, the related eigenvalue problems, warped product manifolds and the Steklov problem on these manifolds. In Section~\ref{sec3} we give the proofs of Theorems~\ref{thm1}, \ref{thm2}, \ref{thm3} and \ref{thm4} and the proof sketch for Remark~\ref{rem6}. In Section~\ref{sec4} we discuss related further problems.

%
%

\section{Preliminaries}\label{sec2}

In this section we collect some preparatory materials which are needed in the proofs of our results, including differential forms, relevant eigenvalue problems on differential forms, warped product manifolds and the characterization of eigenvalues and eigenforms of Steklov problems on warped product manifolds.

\subsection{Differential forms}

Let $(M^{n},g)$ be an $n$-dimensional ($n\geq 2$) smooth compact connected oriented Riemannian manifold with boundary $\Sigma=\pt M$. Let $\nu$ be the unit outward normal vector field along $\Sigma$. Let $J:\Sigma \rightarrow M$ be the inclusion mapping. Let $\Omega^p(M)$ be the space of sections of the vector bundle $\Lambda^p(T^*M)$ of smooth exterior differential forms of degree $p$. Let $d$ and $\delta$ be the differential and co-differential operators on $\Omega^p(M)$, respectively. Let $\mathcal{H}^p(M)$, $\mathcal{E}^p(M)$, $c\mathcal{E}^p(M)$, $\mathcal{C}^p(M)$, and $c\mathcal{C}^p(M)$ denote respectively the subspace in $\Omega^p(M)$ of harmonic fields (i.e., closed and co-closed forms), exact forms, co-exact forms, closed forms, and co-closed forms. Similarly we can define spaces and operators on the boundary $\Sigma=\pt M$, e.g., $\Omega^p(\Sigma)$, $d^\Sigma$, $\delta^\Sigma$. We use $|\cdot|$ to designate the pointwise norm on the differential forms induced by the Riemannian metric $g$, and use $dv$ and $da$ to denote the volume element and area element on $M$ and $\Sigma$, respectively.

Recall that the Hodge Laplacian on differential forms is given by
\begin{equation}
\Delta \omega =d\delta \omega+\delta d\omega,\ \omega \in \Omega^p(M).
\end{equation}
Then we have the classical Bochner formula
\begin{equation}
\Delta \omega =\nabla^*\nabla \omega+W_{M}^{[p]}(\omega),
\end{equation}
where $\nabla^*\nabla$ denotes the connection Laplacian on differential forms and $W_{M}^{[p]}:\Omega^p(M)\to \Omega^p(M)$ is the Weitzenb\"{o}ck curvature operator given by (see e.g. \cite[Thm.~9.4.1]{Pet16})
\begin{align*}
&W_M^{[p]}(\omega)(X_1,X_2,\dots,X_p)\\
=\sum_{i,j}(Rm(e_i,X_j)&(\omega))(X_1,X_2,\dots,X_{j-1},e_i,X_{j+1},\dots, X_p),\ \{X_k\}_{k=1}^p\subset \mathfrak{X}(M).
\end{align*}
Here $Rm(X,Y)=\nabla_X\nabla_Y-\nabla_Y\nabla_X-\nabla_{[X,Y]}:\Omega^p(M)\to \Omega^p(M)$ denotes the induced Riemann curvature operator acting on the differential forms and $\{e_i\}_{i=1}^{n}$ is a local orthonormal frame for the space $\mathfrak{X}(M)$ of smooth vector fields. By Gallot and Meyer's work \cite{GM75}, if the eigenvalues of the curvature operator are bounded below by $\gamma\in \R$, then
\begin{equation*}
\langle W_M^{[p]}(\omega),\omega\rangle \geq p(n-p)\gamma |\omega|^2,\quad \omega \in \Omega^p(M).
\end{equation*}

On the boundary $\Sigma=\pt M$, we use $S:\mathfrak{X}(\Sigma)\to \mathfrak{X}(\Sigma)$ to denote the Weingarten operator (shape operator) with respect to $-\nu$, i.e., $S(X)=\nabla_X \nu$ for $X\in \mathfrak{X}(\Sigma)$. Then the eigenvalues of $S$ are the principal curvatures of the boundary $\Sigma$ in $M$.

\subsection{Eigenvalue problems}\label{sec2.2}
In this subsection we present the Steklov eigenvalue problem for differential forms introduced in \cite{BS08, Kar19}. This eigenvalue problem has some motivation relating to the investigations of Maxwell's equations in mathematical physics. Given $\phi \in \Omega^p(\Sigma)$, there exists a solution $\varphi\in \Omega^p(M)$ satisfying (see \cite[Lemma~3.4.7]{Sch95})
\begin{align}\label{eq-DN}
\begin{cases}
\Delta \varphi=0,\ \delta \varphi  =0\text{ on }M,\\
J^*\varphi =\phi\text{ on }\Sigma;
\end{cases}
\end{align}
the solution is unique up to an element of $\mathcal{H}_D^p(M)$, the space of harmonic fields with Dirichlet boundary condition. For definiteness, we may fix the solution $\varphi$ such that it is orthogonal to $\mathcal{H}_D^p(M)$ (see Proposition~3.11 in \cite{Kar19}). Then the Dirichlet-to-Neumann operator $\Lambda^{[p]}$ is defined by
\begin{equation*}
\Lambda^{[p]}(\phi):=i_\nu d\varphi\in \Omega^p(\Sigma).
\end{equation*}
By \cite{Kar19}, the followings hold true.
\begin{thm}
The operator $\Lambda^{[p]}:\Omega^p(\Sigma)\to \Omega^p(\Sigma)$ is well-defined and self-adjoint. Moreover, it admits the following properties.
\begin{enumerate}
  \item There hold $\mathrm{ker}(\Lambda^{[p]})=J^*\mathcal{H}^p(M)$ and $\mathcal{E}^p(\Sigma)\subset \mathrm{ker}( \Lambda^{[p]})$.
  \item The restriction on the space of co-closed forms $\Lambda^{[p]}:c\mathcal{C}^p(\Sigma)\to c\mathcal{C}^p(\Sigma)$ is an operator with a compact resolvent. The eigenvalues of the restricted operator consist of a sequence
      \begin{align*}
0\leq \sigma_1^{[p]}(\Sigma)\leq \sigma_2^{[p]}(\Sigma)\leq \dots \nearrow \infty,
\end{align*}
with the account of multiplicities.
  \item The eigenvalues satisfy the variational principle
  \begin{align}
\sigma_k^{[p]}(\Sigma)=\inf_{E_k }\sup_{\varphi \in E_k\setminus \{0\}}\frac{\int_M  |d  \varphi |^2 dv}{\int_\Sigma |J^*\varphi |^2 da},
\end{align}
where $E_k$ ranges over $k$-dimensional subspaces in $\Omega^p(M)$ satisfying $J^*E_k\subset c\mathcal{C}^p(\Sigma)$.
\end{enumerate}
\end{thm}
To show the sharpness of our results, we present the spectrum of $\Lambda^{[p]}$ for the unit Euclidean ball $B^{n}(1)$.
\begin{exam}[{\cite{GM75},\cite{IT78},\cite[Prop.~7]{RS14},\cite[Thm.~8.1]{Kar19}}]\label{exam1}
  Let $M=\overline{B^{n}(1)}$ be the unit Euclidean ball. Let $P_{l,p}$ be the space of homogeneous polynomial $p$-forms of degree $l\geq 0$ in $\R^{n}$. Define three subspaces of $P_{l,p}$ as follows:
  \begin{align*}
  H_{l,p}&=\{\omega\in P_{l,p}|\Delta \omega =0,\ \delta \omega=0\},\\
  H'_{l,p}&=\{\omega\in H_{l,p}|d\omega=0\},\\
  H''_{l,p}&=\{\omega \in H_{l,p}|i_\nu\omega=0\text{ on }\SS^{n-1}\}.
  \end{align*}
  (i) Let $1\leq p\leq n-2$. Then $\mathcal{H}^p(\SS^{n-1})=0$ and $\Omega^p(\SS^{n-1})=\mathcal{E}^p(\SS^{n-1})\oplus c\mathcal{E}^p(\SS^{n-1})$. By \cite{IT78}, we know that $\mathcal{E}^p(\SS^{n-1})=\bigoplus_l(J^*H'_{l,p})$, $c\mathcal{E}^p(\SS^{n-1})=\bigoplus_l(J^*H''_{l,p})$, and
  \begin{equation*}
  \delta^{\SS^{n-1}}:J^*H'_{l,p}\to J^*H''_{l+1,p-1}
  \end{equation*}
  is an isomorphism. So $\mathrm{dim} J^* H''_{1,p}=\mathrm{dim} J^* H'_{0,p+1}=C_{n}^{p+1}$. By Theorem~8.1 in \cite{Kar19}, the spaces $J^* H'_{l-1,p}$ and $J^*H''_{l,p}$ for $l\geq 1$ form the eigenbases of $\Lambda^{[p]}$. Precisely, if $\phi\in J^* H'_{l-1,p}$, then $\Lambda^{[p]}\phi=0$; if $\phi\in J^* H''_{l,p}$, then $\Lambda^{[p]}\phi=(p+l)\phi$. In particular, the first positive eigenvalue of $\Lambda^{[p]}$ is $p+1$ with multiplicity $C_{n}^{p+1}$. (ii) Let $p=n-1$. If $\phi\in J^* H'_{l-1,n-1}$, then $\Lambda^{[n-1]}\phi=0$; if $\phi$ is the volume form of $\SS^{n-1}$, then $\Lambda^{[{n-1}]}\phi=n\phi$.
\end{exam}

\begin{rem}\label{rem9}
In Remark~2.6 of \cite{Kar19} (cf. (iii) on Page~112 of \cite{CGGS24}, and Remark~2.1 of \cite{Cha25}), it is stated that $I_{n-1}=1$ (see \eqref{eq-I} for its definition) and $\Lambda^{[n-1]} \equiv 0$ on $\Omega^{n-1}(\pt M)$. Here we present some details on this statement.

Let $(M^n,g)$ be an $n$-dimensional connected oriented Riemannian manifold with boundary. Assume that the boundary $\pt M$ has $k$ connected components. We have the long exact de Rham cohomology sequence of the pair $(M,\pt M)$ as follows
\begin{align*}
\cdots \longrightarrow  H^{n-1}_{dR}(M)\longrightarrow H^{n-1}_{dR}(\pt M)\longrightarrow H^n_{dR}(M,\pt M)\longrightarrow H^n_{dR}(M)\longrightarrow \{0\}.
\end{align*}
Note that from right to left we have
\begin{align*}
&H^n_{dR}(M)\simeq \mathcal{H}^n_N(M)\simeq H_0(M,\pt M; \R)\simeq \{0\},\\
& H^n_{dR}(M,\pt M)\simeq \mathcal{H}^n_D(M)\simeq H_0(M; \R)\simeq \R,\\
&H^{n-1}_{dR}(\pt M)\simeq \mathcal{H}^{n-1}(\pt M)\simeq H_0(\pt M; \R)\simeq \R^k.
\end{align*}
Then the dimension of the kernel of the (restricted) Dirichlet-to-Neumann map $\Lambda^{[n-1]}$ equals
\begin{align}\label{eq-I}
I_{n-1}:=\mathrm{dim}\, \mathrm{im}\{J^*:H_{dR}^{n-1}(M)\to H_{dR}^{n-1}(\pt M)\}=k-1.
\end{align}
Since $c\mathcal{C}^{n-1}(\pt M)=\mathcal{H}^{n-1}(\pt M)$ is $k$-dimensional, we see the Dirichlet-to-Neumann map $\Lambda^{[n-1]} \not\equiv 0 $ on $c\mathcal{C}^{n-1}(\pt M)$. Moreover, we claim that the boundary volume form $dv_{\pt M}$ generates the $1$-dimensional space where $\Lambda^{[n-1]}$ acts nontrivially and is the only eigenform (up to a multiplicative constant). Indeed, by Remark~2.3 in \cite{CHL26}, we have the explicit expression for the Dirichlet-to-Neumann map
\begin{align*}
\Lambda^{[n-1]} \varphi =\frac{\int_{\pt M}\varphi}{|M|}dv_{\pt M},\quad \varphi \in \Omega^{n-1}(\pt M).
\end{align*}
In particular, we get
\begin{align*}
\Lambda^{[n-1]}dv_{\pt M} =\frac{|\pt M|}{|M|}dv_{\pt M},
\end{align*}
which means that $\sigma^{[n-1]}_1=|\pt M|/|M|$ is the only eigenvalue of the operator $\Lambda^{[n-1]}$ on $c\mathcal{C}^{n-1}(\pt M)$ and $dv_{\pt M}$ is the corresponding eigenform.
\end{rem}

For later discussion we also need some facts on the boundary Hodge Laplacian $\Delta_\Sigma$, especially when restricted to the space of co-closed differential forms
\begin{align*}
\Delta_\Sigma^{[p]}:c\mathcal{C}^p(\Sigma)\to c\mathcal{C}^p(\Sigma).
\end{align*}
This restricted operator is a non-negative self-adjoint elliptic operator with eigenvalues
\begin{align*}
0\leq \lambda_1^{[p]}(\Sigma)\leq \lambda_2^{[p]}(\Sigma)\leq \dots \nearrow \infty.
\end{align*}
The eigenvalue $\lambda_k^{[p]}(\Sigma)$ admits the variational principle
\begin{align}
\lambda_k^{[p]}(\Sigma)=\inf_{F_k\subset c\mathcal{C}^p(\Sigma)}\sup_{\phi \in F_k\setminus \{0\}}\frac{\int_\Sigma  |d^\Sigma \phi |^2 da}{\int_\Sigma  |\phi |^2 da},
\end{align}
where $F_k$ ranges over $k$-dimensional subspaces of $c\mathcal{C}^p(\Sigma)$.

For readers' convenience, we describe the spectrum of the operator $\Delta_\Sigma^{[p]}$ on the whole space $\Omega^p(\Sigma)$ for the unit Euclidean ball.
\begin{exam}[{\cite{GM75},\cite{IT78},\cite[Prop.~7]{RS14},\cite[Thm.~8.1]{Kar19}}]\label{exam2}
  Let $M=\overline{B^{n}(1)}$ be the unit Euclidean ball. (i) Let $1\leq p\leq n-2$. The spaces $J^* H'_{l-1,p}$ and $J^*H''_{l,p}$ for $l\geq 1$ form the eigenbases of $\Delta_\Sigma^{[p]}$. Precisely, if $\phi\in J^* H'_{l-1,p}$, then
  \begin{equation*}
  \Delta_\Sigma^{[p]}\phi= (l+p-1)(n+l-p-1) \phi;
  \end{equation*}
  if $\phi\in J^* H''_{l,p}$, then
    \begin{equation}\label{eq-cc}
  \Delta_\Sigma^{[p]}\phi= (l+p)(n+l-p-2) \phi.
  \end{equation}
  (ii) Let $p=n-1$. If $\phi\in J^* H'_{l-1,n-1}$, then
  \begin{equation*}
  \Delta_\Sigma^{[n-1]}\phi= (l+n-2)l \phi;
  \end{equation*}
  if $\phi$ is the volume form of $\SS^{n-1}$, then $\Delta_\Sigma^{[n-1]}\phi=0$.
\end{exam}

\subsection{Warped product manifolds}\label{sec2.3}

Let $M=[0,R]\times\mathbb{S}^{n-1}$ be an $n$-dimensional warped product manifold with the metric
\begin{equation*}
    g=dr^{2}+h^{2}(r)g_{\mathbb{S}^{n-1}}.
\end{equation*}
Assume that the warping function $h(r)$ satisfies:
\begin{align}\label{eq-A}
    h\in C^{\infty}([0,R]);\ h(r)>0,\ \forall r\in(0,R];\ h'(0)=1, h^{(2k)}(0)=0,\ \forall k\geq 0.
\end{align}

Let $\{\partial_{r}, E_{1},\cdots,E_{n-1}\}$ be an orthonormal basis of the tangent space $T_{(r,x)}M$. Then the sectional curvatures at the point $(r,x)$ are given by (see \cite{Pet16})
\begin{align*}
    K(\partial_{r}\wedge E_{i})&=-\frac{h''(r)}{h(r)},\ i=1,\cdots,n-1,\\
    K(E_{i}\wedge E_{j})&=\frac{1-(h'(r))^{2}}{h(r)^{2}},\ i\neq j.
\end{align*}
And the Ricci curvature of the warped product manifold $M$ reads as
\begin{equation*}
    \mathrm{Ric}_{g}=-\left(\dfrac{h''}{h}-(n-2)\dfrac{1-(h')^{2}}{h^{2}}\right)g-(n-2)\left(\dfrac{h''}{h}+\dfrac{1-(h')^{2}}{h^{2}}\right)dr^{2}.
\end{equation*}
On the other hand, the boundary $\partial M=\{R\}\times\mathbb{S}^{n-1}$ of $M$ is umbilical and its principal curvatures are
\begin{equation*}
    \kappa_{1}=\cdots=\kappa_{n-1}=\frac{h'(R)}{h(R)}.
\end{equation*}
To ensure that the warped product manifold $M$ has nonnegative Ricci curvature and a strictly convex boundary, we have the following lemma.
\begin{lem}[\cite{Xio22}]\label{lem-h}
If the warped product manifold $M$ above has nonnegative Ricci curvature and a strictly convex boundary, then the warping function $h(r)$ verifies
	\begin{equation*}
		h''(r)\leq 0, \ 0<h'(r)\leq 1,\ r\in[0,R].
	\end{equation*}
\end{lem}
\begin{rem}\label{rem-curv}
For warped product manifolds with $\mathrm{Ric}_g\geq 0$ and strictly convex boundary, its curvature operator is nonnegative, so all Weitzenb\"{o}ck curvatures are nonnegative. Indeed, using the expressions of the curvature operator on Page 121 of Petersen's book \cite{Pet16}, we see that the curvature operator has only two eigenvalues $-h''/h$ and $(1-(h')^2)/h^2$ (with multiplicities). These two eigenvalues are nonnegative by Lemma~\ref{lem-h}. Thus the curvature operator is nonnegative.
\end{rem}

\subsection{The eigenvalues and eigenforms of the Dirichlet-to-Neumann operator on warped product manifolds}

In \cite{Cha25}, Chakradhar gives the characterization of the Steklov eigenvalues and eigenforms for warped product manifolds in consideration.
\begin{prop}[\cite{Cha25}]\label{prop-ODE}
Let $M=[0,R]\times\mathbb{S}^{n-1}$ be an $n$-dimensional warped product manifold with the metric
\begin{equation*}
    g=dr^{2}+h^{2}(r)g_{\mathbb{S}^{n-1}},
\end{equation*}
where $h(r)$ satisfies the smoothness assumption~\eqref{eq-A}, and $p\in\{1,\cdots,n-2\}$. Then the eigenforms of the operator $\Lambda^{[p]}$ are of the following form of separation of variables
\begin{equation*}
        \varphi(r,x)=\psi(r)\omega(x),
\end{equation*}
where $\omega$ is a co-closed eigenform of the Hodge Laplacian $\Delta^{[p]}_{\SS^{n-1}}$ corresponding to some eigenvalue $\lambda^{[p]}_{(m)}=(m+p)(n+m-p-2)$ with degree $m\in \N$ (see \eqref{eq-cc}), and $\psi(r)$ is a nontrivial solution of the ODE
    \begin{align}\label{radial-ODE}
        \begin{cases}
            \psi''(r)+(n-2p-1)\dfrac{h'(r)}{h(r)}\psi'(r)-\lambda^{[p]}_{(m)}\dfrac{\psi(r)}{h(r)^{2}}=0,\ r\in(0,R],\\
            \psi(0)=0.
        \end{cases}
    \end{align}
Moreover, the Steklov eigenvalue associated with $\varphi$ is given by
    \begin{equation}
        \sigma^{[p]}_{(m)}=\frac{\psi'(R)}{\psi(R)}.
    \end{equation}
\end{prop}
\begin{rem}\label{rem13}
In fact, one may check that Chakradhar's argument also works for the degree $p=n-1$. For this case, the only co-closed eigenform of the Hodge Laplacian on the boundary is the volume form and correspondingly $\lambda^{[n-1]}_{(1)}=0$ (see Example~\ref{exam2}). Therefore $\psi$ satisfies
\begin{align*}
\psi''(r)-(n-1)\dfrac{h'(r)}{h(r)}\psi'(r)=0,\ r\in(0,R],\text{ with }\psi(0)=0.
\end{align*}
Then we solve it to get
\begin{align*}
\psi'(r)=ch^{n-1}(r),\quad \psi(r)=c\int_0^r h^{n-1}(s)ds
\end{align*}
for some $c\neq 0$ and consequently
\begin{align*}
\sigma^{[n-1]}_{(1)}=\sigma^{[n-1]}_1=\frac{h^{n-1}(R)}{\int_0^R h^{n-1}(r)dr}=\frac{|\pt M|}{|M|}.
\end{align*}
This is consistent with Remark~\ref{rem9}. And with this expression at hand, we can readily prove that
\begin{align*}
n\frac{h'(R)}{h(R)} \leq \sigma^{[n-1]}_{(1)}\leq \frac{n}{h(R)}
\end{align*}
under the conditions of Theorem~\ref{thm1}.
\end{rem}

Based on this proposition, Chakradhar obtained a sharp lower bound for the eigenvalues of the Dirichlet-to-Neumann operator on differential forms for certain degrees on warped product manifolds.
\begin{thm}[\cite{Cha25}]\label{thm-Cha}
Let $M=[0,R]\times\mathbb{S}^{n-1}$ be an $n$-dimensional ($n\geq 3$) warped product manifold with the metric
\begin{equation*}
    g=dr^{2}+h^{2}(r)g_{\mathbb{S}^{n-1}},
\end{equation*}
where $h(r)$ satisfies the smoothness assumption~\eqref{eq-A}. Suppose that $M$ has $\mathrm{Ric}_g\geq 0$ and the boundary $\pt M$ is strictly convex. Denote by $\kappa\in\mathbb{R}_{>0}$ the principal curvature of the boundary $\partial M$. Then for any $1\leq p\leq(n-1)/2$, the $m$th nonzero Steklov eigenvalue for $p$-forms without counting multiplicities satisfies
	\begin{equation*}
		\sigma^{[p]}_{(m)}\geq(m+p)\kappa,\ \forall m\in\mathbb{N}.
	\end{equation*}
For a given $m\in\mathbb{N}$, the equality holds if and only if $M$ is a Euclidean ball of radius $1/\kappa$.
\end{thm}

Note that here $\kappa=h'(R)/h(R)$. As mentioned in Section~\ref{sec1}, we extend Chakradhar's result to the full range $p\in [1,n-2]$, and improve on his result for small degrees $p\in [1,(n-2)/2]$.

\section{The proofs of main results}\label{sec3}

In this section, we give the proofs of Theorems~\ref{thm1}, \ref{thm2}, \ref{thm3} and \ref{thm4}.

First we study some properties of the radial function $\psi(r)$ of the eigenform for $\Lambda^{[p]}$ in Proposition~\ref{prop-ODE}.
\begin{prop}\label{prop-psi}
Let $\psi(r)$ be the radial function satisfying the equation~\eqref{radial-ODE} in Proposition~\ref{prop-ODE}. Then up to a sign $\psi(r)$ satisfies
    \begin{equation*}
        \psi(r)>0,\ \psi'(r)>0,\ \forall  r\in (0,R].
    \end{equation*}
    Moreover, we have $\psi^{(k)}(0)=0$ for $k=1,2,\dots, m+p-1$; in other words, $\psi(r)=O(r^{m+p})$ as $r\to 0$.
\end{prop}
\begin{rem}
Throughout the paper, we assume that $\psi$ satisfies the properties in Proposition~\ref{prop-psi} unless otherwise stated.
\end{rem}
\begin{proof}
First, since $\psi(r)$ is a nontrivial solution of the equation~\eqref{radial-ODE}, we may assume $\psi(r_0)>0$ for some point $r_0\in (0,R]$. Then by the maximum principle, $\psi(R)=\max \psi>0$ and $\psi(0)=\min \psi=0$. Moreover, by the strong maximum principle, we see $\psi(r)>0$ for $r\in (0,R]$.

Second, fix any $r\in (0, R]$. By Lagrange Mean Value Theorem, there exists a point $r_1\in (0,r)$ such that $\psi'(r_1)>0$. Then integrating the equation
\begin{equation*}
\dfrac{d}{dr}\left(h^{n-2p-1}(r)\psi'(r)\right)=\lambda^{[p]}_{(m)}\psi(r)h^{n-2p-3}(r)
\end{equation*}
over $[r_1,r]$, we see $\psi'(r)>0$.

Third, recall that $\psi$ satisfies
\begin{equation*}
h\psi''=\lambda^{[p]}_{(m)}\frac{\psi}{h}-(n-2p-1)h'\psi'.
\end{equation*}
Letting $r\to 0$ and applying L'H\^{o}pital's rule, we get
\begin{equation*}
0=\lambda^{[p]}_{(m)}\psi'(0)-(n-2p-1)\psi'(0),
\end{equation*}
which implies $\psi'(0)=0$. So if $m+p=2$, we finish the proof. In the following we consider the case $m+p\geq 3$.

Since $\psi(r)$ is a smooth function on $[0,R]$, we may assume
\begin{align*}
\psi(r)&=\frac{\psi^{(k)}(0)}{k!}r^k+o(r^k),\\
\psi'(r)&=\frac{\psi^{(k)}(0)}{(k-1)!}r^{k-1}+o(r^{k-1}),\\
\psi''(r)&=\frac{\psi^{(k)}(0)}{(k-2)!}r^{k-2}+o(r^{k-2}),
\end{align*}
as $r\to 0$, where $k\geq 2$ is an integer. We shall prove that for $k=2,\dots, m+p-1$, $\psi^{(k)}(0)=0$.

In fact, using the expansions for $\psi(r)$, $\psi'(r)$, $\psi''(r)$ above, $h(r)=r+o(r)$ and $h'(r)=1+o(1)$ as $r\to 0$ in the equation
\begin{equation*}
h^2\psi''=\lambda^{[p]}_{(m)} \psi-(n-2p-1)hh'\psi',
\end{equation*}
we may derive
\begin{align*}
\frac{\psi^{(k)}(0)}{(k-2)!}r^{k}=\lambda^{[p]}_{(m)}\frac{\psi^{(k)}(0)}{k!}r^k-(n-2p-1)\frac{\psi^{(k)}(0)}{(k-1)!}r^k+o(r^{k}),\text{ as }r\to 0.
\end{align*}
So we get
\begin{align*}
\psi^{(k)}(0)\left(k(k-1)+(n-2p-1)k-\lambda^{[p]}_{(m)}\right)=0,
\end{align*}
or
\begin{align*}
\psi^{(k)}(0)(k-m-p)(k+n+m-2-p)=0.
\end{align*}
Consequently, letting $k=2,\dots, m+p-1$ in turn, we find $\psi^{(k)}(0)=0$ for these $k$ and finish the proof.
\end{proof}

\subsection{The proof of Theorem~\ref{thm1}}

\begin{proof}
By Proposition~\ref{prop-ODE}, the eigenvalue $\sigma^{[p]}_{(m)}$ is given by
\begin{equation*}
\sigma^{[p]}_{(m)}=\frac{\psi'_{m}(R)}{\psi_{m}(R)},
\end{equation*}
where $\psi_m(r)$ is a nontrivial solution to the ODE~\eqref{radial-ODE}. Let
\begin{equation*}
q(r):=h(r)\psi_{m}'(r)-(m+p)h'(r)\psi_{m}(r).
\end{equation*}
By Proposition~\ref{prop-psi}, $q(0)=0$. Taking derivative on $q$ and using the ODE
\begin{equation*}
\psi''_{m}(r)=\lambda^{[p]}_{(m)}\frac{\psi_{m}(r)}{h(r)^{2}}-(n-2p-1)\frac{h'(r)}{h(r)}\psi'_{m}(r),
\end{equation*}
we obtain
\begin{align}
        q'=&h'\psi_{m}'+h\psi_{m}''-(m+p)h''\psi_m-(m+p)h'\psi_{m}' \nonumber \\
        \geq &(1-m-p)h'\psi_{m}'+\left(\lambda^{[p]}_{(m)}\frac{\psi_{m}}{h}-(n-2p-1)h'\psi'_{m}\right)\nonumber \\
        =&(2+p-m-n)h'\psi_{m}'+\lambda^{[p]}_{(m)}\frac{\psi_{m}}{h}\nonumber \\
        =&\lambda^{[p]}_{(m)}\left(\frac{\psi_{m}}{h}-\frac{1}{m+p}h'\psi_{m}'\right), \label{eq-q'}
\end{align}
where we used $h''\leq 0$. Define
\begin{align*}
u(r):=(m+p)\frac{\psi_{m}(r)}{h(r)}-h'(r)\psi_{m}'(r).
\end{align*}
Then we see
\begin{align*}
u(0)=\left(m+p-1\right)\psi_{m}'(0)=0.
\end{align*}
Furthermore, we derive
\begin{align*}
u'&=(m+p)\frac{\psi'}{h}-(m+p)\frac{\psi}{h^2}h'-h''\psi'-h'\psi''\\
&\geq (m+p)\frac{\psi'}{h}-(m+p)\frac{\psi}{h^2}h'\\
&\quad -h'\left(\lambda^{[p]}_{(m)}\frac{\psi_{m}}{h^{2}}-(n-2p-1)\frac{h'}{h}\psi'_{m}\right)\\
&=\left(m+p+(n-2p-1)(h')^2\right)\frac{\psi'}{h}-(m+p+\lambda^{[p]}_{(m)})\frac{h'\psi}{h^2}.
\end{align*}
Then by the definition of $u$, we get
\begin{align*}
u'&\geq \left(m+p+(n-2p-1)(h')^2\right)\frac{1}{hh'}\left((m+p)\frac{\psi_{m}}{h}-u\right)\\
&\quad -(m+p+\lambda^{[p]}_{(m)})\frac{h'\psi}{h^2}\\
&=-\left(m+p+(n-2p-1)(h')^2\right)\frac{1}{hh'}u+\frac{\psi}{h^2h'}(m+p)^2(1-(h')^2)\\
&\geq -\left(m+p+(n-2p-1)(h')^2\right)\frac{1}{hh'}u.
\end{align*}
Define
\begin{align*}
c(r):= \left(m+p+(n-2p-1)(h'(r))^2\right)\frac{1}{h(r)h'(r)}.
\end{align*}
We find
\begin{align*}
c(r)=(n+m-p-1)\frac{1}{r}(1+o(1)),\text{ as } r\to 0,
\end{align*}
and
\begin{align*}
u'(r)+c(r)u(r)\geq 0.
\end{align*}
Then for $\varepsilon\in (0,R]$, we may get
\begin{align*}
\frac{d}{dr}\left(e^{\int_\varepsilon^r c(s)ds}u(r)\right)\geq 0,
\end{align*}
which after integration over $[\varepsilon,r]$ implies
\begin{align*}
e^{\int_\varepsilon^r c(s)ds}u(r)\geq u(\varepsilon).
\end{align*}
As a result, we obtain
\begin{align*}
u(r)\geq e^{-\int_\varepsilon^r c(s)ds}u(\varepsilon)\to 0,\text{ as } \varepsilon\to 0,
\end{align*}
since $\int_\varepsilon^r c(s)ds\to +\infty$ and $u(\varepsilon)\to 0$ as $\varepsilon \to 0$. Then $q'\geq \lambda^{[p]}_{(m)}u/(m+p)\geq 0$ and $q\geq 0$. So we finish the proof of the inequality.

When $M$ is a Euclidean ball, we know $h(r)=r$, and by Theorem~8.1 in \cite{Kar19}, $\sigma^{[p]}_{(m)}=(m+p)/R$. So the equality holds.

Conversely, when the equality holds, i.e., $q(R)=0$, we get $q\equiv 0$ and $q'\equiv 0$. Then the inequality in \eqref{eq-q'} becomes an equality, forcing $h''\equiv 0$. Thus we conclude $h(r)=r$, i.e., $M$ is isometric to a Euclidean ball of radius $R$. So we complete the proof.
\end{proof}

\subsection{The proof of Theorem~\ref{thm2}}

\begin{proof}
By Proposition~\ref{prop-ODE}, the eigenvalue $\sigma^{[p]}_{(m)}$ is given by
\begin{equation*}
\sigma^{[p]}_{(m)}=\frac{\psi'_{m}(R)}{\psi_{m}(R)},
\end{equation*}
where $\psi_m(r)$ is a nontrivial solution to the ODE~\eqref{radial-ODE}. Let
\begin{equation*}
q(r):=h(r)\psi_{m}'(r)-(m+p)\psi_{m}(r).
\end{equation*}
By Proposition~\ref{prop-psi}, $q(0)=0$. Taking derivative on $q$ and using the ODE
\begin{equation*}
\psi''_{m}(r)=\lambda^{[p]}_{(m)}\frac{\psi_{m}(r)}{h(r)^{2}}-(n-2p-1)\frac{h'(r)}{h(r)}\psi'_{m}(r),
\end{equation*}
we obtain
\begin{align*}
        q'=&h'\psi_{m}'+h\psi_{m}''-(m+p)\psi_{m}'\\
        =&h'\psi_{m}'+h\left(\frac{\lambda^{[p]}_{(m)}}{h^{2}}\psi_{m}-(n-2p-1)\frac{h'}{h}\psi_{m}'\right)-(m+p)\psi_{m}'\\
        =&\frac{\lambda^{[p]}_{(m)}}{h}\psi_{m}-(n-2p-2)h'\psi_{m}'-(m+p)\psi_{m}'.
\end{align*}
By the definition of $q$, we get
    \begin{equation*}
        \psi_{m}'=\frac{1}{h}\left(q+(m+p)\psi_{m}\right).
    \end{equation*}
Hence
    \begin{align*}
        q'=&\frac{\lambda^{[p]}_{(m)}}{h}\psi_{m}-(n-2p-2)\frac{h'}{h}\left(q+(m+p)\psi_{m}\right)-(m+p)\frac{1}{h}\left(q+(m+p)\psi_{m}\right)\\
        =&-\left(\frac{m+p}{h}+(n-2p-2)\frac{h'}{h}\right)q\\
        &+\left(\lambda^{[p]}_{(m)}-(m+p)(n-2p-2)h'-(m+p)^{2}\right)\frac{\psi_{m}}{h}\\
=&-\left(\frac{m+p}{h}+(n-2p-2)\frac{h'}{h}\right)q+(m+p)\left(n-2p-2\right)(1-h')\frac{\psi_{m}}{h}\\
        \geq& -\left(\frac{m+p}{h}+(n-2p-2)\frac{h'}{h}\right)q,
    \end{align*}
where we used $0< h'\leq 1$ and $1\leq p\leq (n-2)/2$. So we have $q\geq 0$ by a similar argument as in the proof of Theorem~\ref{thm1} and get the inequality.

When $M$ is a Euclidean ball, we know $h(r)=r$, and by Theorem~8.1 in \cite{Kar19}, $\sigma^{[p]}_{(m)}=(m+p)/R$. So the equality holds.

Conversely, for $p<(n-2)/2$, when the equality holds, the inequality in the above argument becomes an equality. Then we must have $h'(r)\equiv 1$ and $h(r)=r$, i.e., $M$ is isometric to a Euclidean ball of radius $R$.

Last, if $p=(n-2)/2$ is an integer, then $q'=-(m+p)q/h$. Once again, by a similar argument as in the proof of Theorem~\ref{thm1}, we can see $q\equiv 0$. By Proposition~\ref{prop-ODE}, this gives
\begin{equation*}
        \sigma^{[p]}_{(m)}=\frac{\psi'_{m}(R)}{\psi_{m}(R)}=\frac{m+p}{h(R)}, \ \forall m\geq 1.
\end{equation*}
So we finish the proof.
\end{proof}

\subsection{The proof of Theorem~\ref{thm3}}

\begin{proof}
Let $\psi_{m+1}$ and $\psi_{m}\in C^{\infty}([0,R])$ be the radial function of the eigenform in Proposition~\ref{prop-ODE} corresponding to the eigenvalue $\sigma^{[p]}_{(m+1)}$ and $\sigma^{[p]}_{(m)}$, respectively. Set
    \begin{equation}
        q(r):=h(r)\psi_{m+1}'(r)-h(r)\frac{\psi_{m}'(r)}{\psi_{m}(r)}\psi_{m+1}(r)-\psi_{m+1}(r).
    \end{equation}
First we claim that
\begin{align}\label{eq-claim}
\lim_{r \to 0} h(r)\frac{\psi_{m}'(r)}{\psi_{m}(r)}=m+p.
\end{align}
In fact, note that $\psi_m'(r)/\psi_m(r)$ is the $m$th Steklov eigenvalue without counting multiplicities on the warped product manifold of radius $r$. When $r$ is small enough, the warped product manifold of radius $r$ is approximately a Euclidean ball of radius $r$ which has the $m$th Steklov eigenvalue $(m+p)/r$. So we get
\begin{align*}
\frac{\psi_m'(r)}{\psi_m(r)}=\frac{m+p}{r}(1+o(1)),\text{ as }r\to 0.
\end{align*}
Then we can readily get the claim. (Alternatively, one may prove the claim by using L'H\^{o}pital's rule and the ODE satisfied by $\psi_m$.)

Now by Proposition~\ref{prop-psi} we know $q(0)=0$. Taking the derivative on $q$, and using the equations
    \begin{align*}
        \psi''_{m}(r)&=\lambda^{[p]}_{(m)}\frac{\psi_{m}(r)}{h(r)^{2}}-(n-2p-1)\frac{h'(r)}{h(r)}\psi'_{m}(r),\\
         \psi''_{m+1}(r)&=\lambda^{[p]}_{(m+1)}\frac{\psi_{m+1}(r)}{h(r)^{2}}-(n-2p-1)\frac{h'(r)}{h(r)}\psi'_{m+1}(r),
    \end{align*}
we may deduce
    \begin{align*}
        q'=&(h'-1)\psi_{m+1}'+h\psi_{m+1}''-h'\frac{\psi_{m}'}{\psi_{m}}\psi_{m+1}-h\frac{\psi_{m}''}{\psi_{m}}\psi_{m+1}+h\frac{(\psi_{m}')^{2}}{(\psi_{m})^{2}}\psi_{m+1}\\
        &-h\frac{\psi_{m}'}{\psi_{m}}\psi_{m+1}'\\
        =&(h'-1)\psi_{m+1}'+\frac{\lambda^{[p]}_{(m+1)}}{h}\psi_{m+1}-(n-2p-1)h'\psi_{m+1}'-h'\frac{\psi_{m}'}{\psi_{m}}\psi_{m+1}\\
        &-\frac{\lambda^{[p]}_{(m)}}{h}\psi_{m+1}+(n-2p-1)h'\frac{\psi_{m}'}{\psi_{m}}\psi_{m+1}+h\frac{(\psi_{m}')^{2}}{(\psi_{m})^{2}}\psi_{m+1}-h\frac{\psi_{m}'}{\psi_{m}}\psi_{m+1}'\\
        =&-((n-2p-2)h'+1)\psi_{m+1}'+\frac{1}{h}\left(\lambda^{[p]}_{(m+1)}-\lambda^{[p]}_{(m)}\right)\psi_{m+1}\\
        &+(n-2p-2)h'\frac{\psi_{m}'}{\psi_{m}}\psi_{m+1}
        +h\frac{(\psi_{m}')^{2}}{(\psi_{m})^{2}}\psi_{m+1}-h\frac{\psi_{m}'}{\psi_{m}}\psi_{m+1}'.
    \end{align*}
By the definition of the function $q$, we get
    \begin{equation*}
        \psi_{m+1}'=\frac{1}{h}\left(q+h\frac{\psi_{m}'}{\psi_{m}}\psi_{m+1}+\psi_{m+1}\right).
    \end{equation*}
Plugging it into the expression for $q'$ we obtain
    \begin{align}
        q'=&-\frac{1}{h}\left(1+(n-2p-2)h'+h\frac{\psi_{m}'}{\psi_{m}}\right)q\nonumber\\
        &+\left((n-2p-2)(1-h')+2m+2p-2h\frac{\psi_{m}'}{\psi_{m}}\right)\frac{\psi_{m+1}}{h}.
    \end{align}
Now define
    \begin{equation}
        u(r):=\left((n-2p-2)(1-h'(r))+2m+2p\right)\psi_{m}(r)-2h(r)\psi_{m}'(r).
    \end{equation}
First we see $u(0)=0$. Taking derivative on $u$ we get
    \begin{align*}
        u'=&-(n-2p-2)h''\psi_{m}+((n-2p-2)(1-h')+2m+2p)\psi_{m}'\\
        &-2h'\psi_{m}'-2h\psi_{m}''\\
        \geq&((n-2p-2)(1-h')+2m+2p)\psi_{m}'-2h'\psi_{m}'-2h\psi_{m}''\\
        =&((n-2p-2)(1-h')+2m+2p)\psi_{m}'-2h'\psi_{m}'\\
        &-2h\left(\frac{\lambda^{[p]}_{(m)}}{h^2}\psi_{m}-(n-2p-1)\frac{h'}{h}\psi_{m}'\right)\\
        =&-\frac{1}{2h}((n-2p-2)(1-h')+2m+2p)u+\frac{\psi_m}{2h}(n-2p-2)^{2}(1-(h')^{2})\\
        \geq&-\frac{1}{2h}((n-2p-2)(1-h')+2m+2p)u,
    \end{align*}
where the first inequality is due to $h''\leq 0$ and $p\leq (n-2)/2$, while the second one due to $0<h'\leq 1$. So we may argue as before to conclude $u\geq 0$ on $[0,R]$, leading to further $q\geq 0$ on $[0,R]$. Thus we get the inequality.

When $M$ is a Euclidean ball, we know $h(r)=r$, and by Theorem~8.1 in \cite{Kar19}, $\sigma^{[p]}_{(m)}=(m+p)/R$. So the equality holds.

Conversely, for $p<(n-2)/2$, when the equality holds, all the inequalities in the above argument become equalities. Then we must have $h'(r)\equiv 1$ and $h(r)=r$, i.e., $M$ is isometric to a Euclidean ball of radius $R$.

Last, if $p=(n-2)/2$ is an integer, then by Theorem~\ref{thm2}, we have
\begin{equation*}
        \sigma^{[p]}_{(m)}=\frac{m+p}{h(R)}, \ \forall m\geq 1.
\end{equation*}
So we finish the proof.
\end{proof}

\subsection{The proof of Theorem~\ref{thm4}}
The proof will rely on the sharp lower bound for $\sigma^{[p]}_{(m)}$ from Theorem~\ref{thm2}.
\begin{proof}
Let
\begin{align*}
\alpha:=\frac{m+p+1}{m+p}.
\end{align*}
Then we need to prove
    \begin{equation*}
        \frac{\sigma^{[p]}_{(m+1)}}{\sigma^{[p]}_{(m)}}=\frac{\psi_{m+1}'(R)\psi_{m}(R)}{\psi_{m+1}(R)\psi'_{m}(R)}\leq\alpha.
    \end{equation*}
Let
    \begin{equation}
        q(r):=\alpha\psi_{m}'(r)\psi_{m+1}(r)-\psi_{m}(r)\psi_{m+1}'(r).
    \end{equation}
Then first we see $q(0)=0$. Taking the derivative on $q$, and using the ODEs satisfied by $\psi_{m}$ and $\psi_{m+1}$ respectively, we derive
    \begin{align*}
        q'=&\alpha\psi_{m}''\psi_{m+1}+\alpha\psi_{m}'\psi_{m+1}'-\psi_{m}'\psi'_{m+1}-\psi_{m}\psi_{m+1}''\\
        =&(\alpha-1)\psi_{m}'\psi_{m+1}'+(\alpha\lambda^{[p]}_{(m)}-\lambda^{[p]}_{(m+1)})\frac{\psi_{m}\psi_{m+1}}{h^{2}}-(n-2p-1)\frac{h'}{h}q\\
        =&(\alpha-1)\psi_{m}'\frac{1}{\psi_{m}}(\alpha\psi_{m}'\psi_{m+1}-q)+(\alpha\lambda^{[p]}_{(m)}-\lambda^{[p]}_{(m+1)})\frac{\psi_{m}\psi_{m+1}}{h^{2}}\\
        &-(n-2p-1)\frac{h'}{h}q      \\
        =&-\left((n-2p-1)\frac{h'}{h}+(\alpha-1)\frac{\psi_{m}'}{\psi_{m}}\right)q\\
        &+\left(\alpha(\alpha-1)\frac{(\psi_{m}')^{2}}{(\psi_{m})^{2}}+\frac{\alpha\lambda^{[p]}_{(m)}-\lambda^{[p]}_{(m+1)}}{h^{2}}\right)\psi_{m}\psi_{m+1},
    \end{align*}
where in the third equality we used the definition of the function $q$. Note that $\psi_m'(r)/\psi_m(r)$ is the $m$th Steklov eigenvalue without counting multiplicities on the warped product manifold of radius $r$. So by Theorem~\ref{thm2}, we find
    \begin{align*}
       \frac{(\psi_{m}')^{2}}{(\psi_{m})^{2}}\geq\frac{(m+p)^{2}}{h^{2}}
       =-\frac{1}{\alpha(\alpha-1)}\frac{\alpha\lambda^{[p]}_{(m)}-\lambda^{[p]}_{(m+1)}}{h^{2}}.
    \end{align*}
Therefore we get
\begin{align*}
        q'&\geq -\left((n-2p-1)\frac{h'}{h}+(\alpha-1)\frac{\psi_{m}'}{\psi_{m}}\right)q.
\end{align*}
Note again that $\psi_m'(r)/\psi_m(r)$ is the $m$th Steklov eigenvalue without counting multiplicities on the warped product manifold of radius $r$, which is equivalent to $(m+p)/r$ as $r\to 0$. As a result, we can deduce
\begin{align*}
(n-2p-1)\frac{h'}{h}+(\alpha-1)\frac{\psi_{m}'}{\psi_{m}}=\frac{n-2p}{r}(1+o(1)),\text{ as } r\to 0.
\end{align*}
Then we may argue as before to derive $q\geq 0$. So we get the inequality.

When $M$ is a Euclidean ball, we know $h(r)=r$, and by Theorem~8.1 in \cite{Kar19}, $\sigma^{[p]}_{(m)}=(m+p)/R$. So the equality holds.

Conversely, for $p<(n-2)/2$, when the equality holds, we must have
\begin{align*}
       \frac{\psi_{m}'(R)}{\psi_{m}(R)}=\frac{m+p}{h(R)}.
    \end{align*}
Then by Theorem~\ref{thm2}, the manifold $M$ is isometric to a Euclidean ball of radius $R$.

Last, if $p=(n-2)/2$ is an integer, then by Theorem~\ref{thm2}, we have
\begin{equation*}
        \sigma^{[p]}_{(m)}=\frac{m+p}{h(R)}, \ \forall m\geq 1.
\end{equation*}
So we finish the proof.
\end{proof}

\subsection{The proof sketch for Remark~\ref{rem6}}
When $\mathrm{Ric}_{g}\leq 0$, the inequalities in Lemma~\ref{lem-h} are reversed. So we have
    \begin{equation*}
        0\leq h''(r),\  1\leq h'(r) ,\ \forall r\in(0,R].
    \end{equation*}
Then all the inequalities in the proofs of Theorems~\ref{thm1}, \ref{thm2}, \ref{thm3} and \ref{thm4} are reversed. Consequently we can get the corresponding conclusions.

\section{Final remarks}\label{sec4}

By Remark~\ref{rem-curv}, for warped product manifolds satisfying the smoothness assumption~\eqref{eq-A} with convex boundary, $\mathrm{Ric}_g \geq 0$ is equivalent to the very strong condition that the curvature operator is nonnegative. Thus it would be an interesting question which results in Theorems~\ref{thm1}, \ref{thm2}, \ref{thm3} and \ref{thm4} still hold under certain weaker curvature conditions for general compact Riemannian manifolds with boundary. For example, one may consider Riemannian manifolds with nonnegative Weitzenb\"{o}ck curvature of some degree. One may also consider weaker convexity of the boundary, e.g., $p$-convexity (see \cite{Xio26} for the definition). For a start, one may first study these problems on Euclidean domains.

Another question is what happens in the setting of Theorem~\ref{thm2} for the remaining degree $(n-2)/2<p\leq n-2$. Note by Remark~\ref{rem2} we have $\sigma^{[n-1]}_{(1)}\leq n/h(R)$ for the highest degree $n-1$. So can we prove
    \begin{equation}
        \sigma^{[p]}_{(m)}\leq \frac{m+p}{h(R)},\ \forall m\geq 1
    \end{equation}
for $(n-2)/2<p\leq n-2$?


\section*{Statements and Declarations}

\noindent \textbf{Funding}: This research was supported by National Key R and D Program of China 2021YFA1001800 and NSFC (Grant no. 12171334).


\noindent \textbf{Conflict of interest}: The authors have no competing interests to declare that are relevant to the content of this article.

\noindent \textbf{Data availability statement}: Data sharing is not applicable to this article as no datasets were generated or
analyzed during the current study.

\noindent \textbf{Authors' contributions}: Both authors contributed to the manuscript equally. Both authors read and approved the final manuscript.

\noindent \textbf{AI Use Disclosure}: During the preparation of this manuscript, AI tools were used to help authors with understanding of known results in literature and English writing.


\bibliographystyle{Plain}

\end{document}